# Switching Dynamics of Shallow Arches

**Priyabrata Maharana**

Department of Mechanical Engineering, Indian Institute of Science, CV Raman Road, Bengaluru, 560012, Karnataka, India.
priyabratam@iisc.ac.in

**Prof. G. K. Ananthasuresh**

Department of Mechanical Engineering, Indian Institute of Science, CV Raman Road, Bengaluru, 560012, Karnataka, India.
suresh@iisc.ac.in

**ABSTRACT**

*This paper presents an analytical method to predict the delayed switching dynamics of nonlinear shallow arches while switching from one state to other state for different loading cases. We study an elastic arch subject to static loading and time-dependent loading separately. In particular, we consider a time-dependent loading that evolves linearly with time at a constant rate. In both the cases, we observed that the switching does not occur abruptly when the load exceeds the static switching load, rather the time scale of the dynamics drastically slows down, hence there is a delay in switching. For time-independent loading this delay increases as the applied load approaches the static switching load. Whereas, for a time-dependent loading the delay is proportional to the rate of the applied load. Other than the loading parameters, the delay switching time also depends on the local curvature of the force-displacement function at the static switching point and the damping co-efficient of the arch material. The delay switching occurs due to the flatness of the energy curve at static switching load. Therefore, we linearize the arch near to the static switching point and get a reduced nonlinear ordinary differential equation to study the switching dynamics of the arch. This reduced equation allows us to derive analytical expressions for the delay switching time of the. We further compare the derived analytical results with the numerical solutions and observed a good agreement between them. Finally, the derived analytical formulae can be used to design arches for a self-offloading dynamic footwear for the diabetics.*

1. Introduction

Arches are the mechanical structures that shows nonlinear force-displacement characteristics when a transverse load is applie. This nonlinear behavior is due to the mutual effect of bending and axial thrust due to the arch deformation. Based on the force-displacement characteristics, arches are of two kinds: bistable and snap-through. For an bistable arch, the force-displacement curve has three zero-force configurations, out of which two are stable states and one is unstable state. On the other hand, the arch is called snap-through when it has only one zero-force state corresponds to the as-fabricated configuration. In both the cases, the arch shows spring softening effect initially, becomes unstable by showing negative stiffness region, and switches to its inverted configuration when the force exceeds a certain threshold value. In past few decades, the use of arches has been increased drastically in fields like micro-relays, electromagnetic actuators, micro-valves, mechanical memory components, retractable devices, consumer products, circuit breakers, and easy-chairs [1-7] due to its unique force-displacement characteristic. Another application that we are interested in is designing a self-offloading footwear for the diabetics using snap-through and bistable arches [8]. The most important factor in designing this kind of footwear is to customize them according to the person's weight and the gait cycle. This customization involves the dynamic analysis of the arches for a certain switching force, and switching and switchback time. Designing of arches taking a pulse load and the subsequent switchback time from the deformed shape to the original undeformed state has already been mentioned in our previous article [wearable]. In this paper, we have covered the switching time of the arch when the load exceeds the threshold value and the behavior of the arch when the load is time varying. The force-time relation of human walking is characterized by a double-hump curve [9], in which the maximum force occurs during flat foot and heel off phase (see Fig. 1) and the lowest during the mid-stance. Though this a nonlinear curve for a complete gait cycle (black-solid curve), every individual gate phase can be approximated as a linear curve of slope that depends on the person's walking speed (magneta dot curve). So, doing these dynamic analyses are important in deciding the arch dimensions more precisely to customize the footwear according to person's weight and gait cycle speed. We show that for a time-dependent load, the critical switching load changes and the switching of the arch delay significantly compared to the static load case.

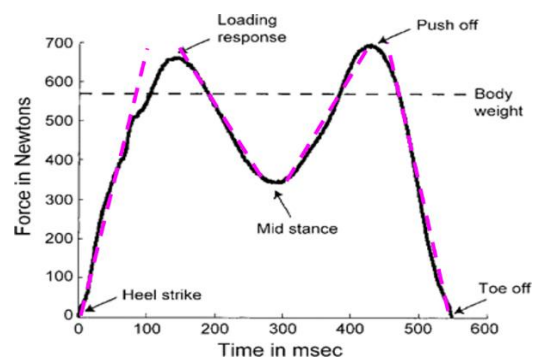

**Fig.1.** The vertical force component of the ground reaction forces during gait.

In past few decades several researchers had found out the switching load by considering it to be quasi-static [10-12], a constant function in time [13] or as an impulse function of time. In their analysis, they assume the arch switches abruptly to its inverted shape when the load slightly higher than the threshold value. Despite of applications of bistable arches for different precission engineering applications, the dynamic behavior is not well studied, specially the dynamics at the time of switching. Therefore, in this paper we analize the local behavior of the arch at the static switching point for different loads, and show that switching of the arch from one stable path to other slows down near the static switching point. Through analytical modeling we show this slowness of switching of arches and verify them using numerical methods. When the applied load is time dependent and vary linearly with time, the switching point of the arch delays by taking more load than the static switching load. This delaying nature of the arch depends on the rate of change of the applied load and the arch geometry. This means, the arch does not switch at the static switching load no matter how slowly we vary the external applied load.

This delay switching occurs at the saddle-node (switching point) due to the nonlinear character of the arch and can also be observed in fields such as, biological science [14-16], statistical mechanics [17], chemical science [18], atmospheric science [19], economics [20, 21], optics [22, 23], ecology [24], in mechanics by end shortening oftlat beam [25], and so on. This phenomena is called the critical slowing down (CSD) [26], and first coined by Racz [17] where he considered the difference between the critical slowing down of linear and nonlinear kinetic Ising model. Bonifacio and Meystre [23], taking the advantage of CSD, thought of novel applications to delay compact optical lines using optical bistability. By taking the controlling parameter independent of time, they showed that the transition of stability slows down drastically when the controlling parameter just above the threshold value and this slowness depends on the closeness of the control parameter to the critical value. Also, when the controlling parameter is a function of time, particularly a ramp function, the switching point of the system gets delayed before it finally switches to another stable point. The point to which the system jumps is called the "tipping point" [25]. This is the point after which the system change its stability by causing a catastrophic failure. Many researchers use this critical transition of the system as an early signal of system failure [25, 27] and can be avoided by the reversal of the control parameter before it reaches tipping point. The time a system takes to reach the tipping point is called the point of no return [28].

In this paper, we focus on the delaying nature of nonlinear shallow arches in two different cases, when force is quasistatic; and in another it is ramp function of time. Through analytical framework, we show that the timescale of the arch decreases compared to the elastic time scale by obeying a power law relation, and the switching delay as the arch approaches the static switching point. We discuss the CSD in detail in Section 2, followed by the theoretical background for the dynamic analysis of the arch in Section 3. In Section 4 we develop a mathematical framework for the switching dynamics, get the analytical formulation for the switching time of the arch for time-independent and time-dependent loading. We extend the method for higher mode approximation in Section 5. The results from both the analytical and numerical methods are presented and discussed in Section 6 followed by a final summary in Section 7.

## 2. Understanding critical slowing down

Due to the nonlinear behavior of arch, there exist multiple displacement values for a paricular force value. When the applied force is less than a threshold force value, $F_c$, there coexist three equilibrium configurations for a force value out of which one is unstable and other two are stable configurations. This leads to one unstable branch and two stable branches on the force-displacement equilibrium plane. The stable equilibrium points corresponds to the minima on the energy-displacement landscape, whereas the maixima corresponds to a unstable equilibrium configurations As the load increases, the stable and unstable point approach each other, coincide, and exchangeing stability at the switching load load $F_c$ due to the saddle node bifurcation. As the force approaches $F_c$, the curvature of the energy-displacement landscape becomes flatter gradually, and finally gets zero at the switching point. When the arch moves through this platue region, its time scale slows down drastically. So, when $F > F_c$ the switching of the arch delay and the arch takes longer time to switch from one stable point to the other. This phenomena is called the critical slowing down of the arch. This occurs due to the zero-stiffness value of the arch at the switching point. It results a higher relaxation time that slows down or delay the dynamics compared to the elastic time scale of the arch. This slowness increases as the arch approaches the static switching point. When the force is time-dependent, this delay in the switching of the arch depends on the rate at which the applied force changes with time. No matter how slowly we vary the load this delay phenomena occurs at the static switching point of the arch.

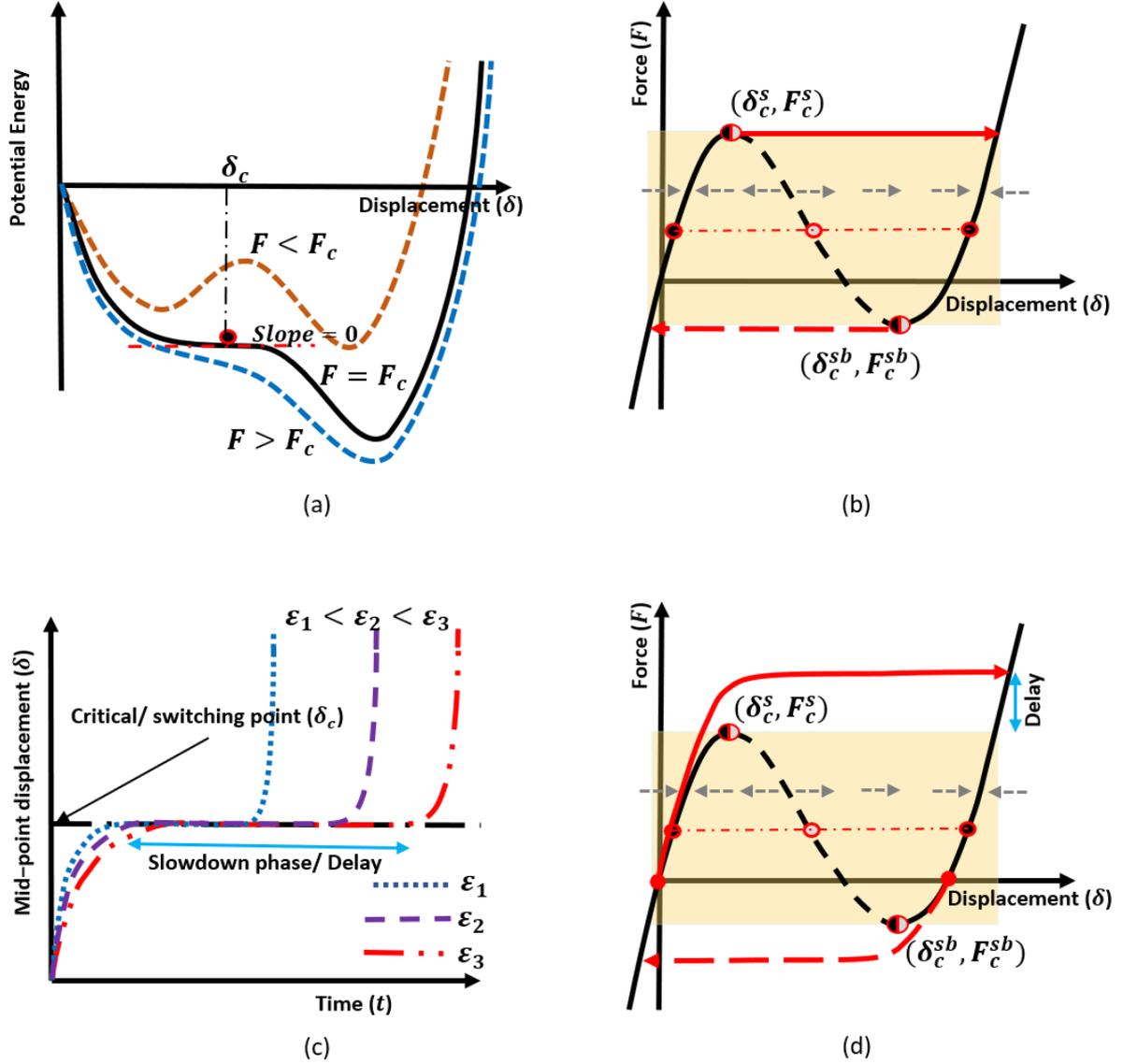

**Fig. 2.** (a) Potential energy-displacement plot for an bistable arch for different loading at the mid-span length of the arch, $F < F_c$ (orange-dash curve), $F = F_c$ (black-solid curve), and $F > F_c$ (blue-dash curve). The black ball indicates the state of the arch, (b) static force-displacement behavior of the arch when a time-independent load is applied with red arrow indicating the switching of the arch after the static switching point, (c) slowdown of arch while switching for different load value when time-independent load is applied. Here $\varepsilon$ indicates the closeness of the external applied load to the static switching force, i.e., $\varepsilon = |F - F_c|$, (d) behavor of the arch when a time-dependent load is applied. The red curve indicates the path of the dynamic switching curve w.r.t. the static force-displacement curve (black curve), and the blue colour arrow indicates the delay in the switching for particular load rate. The yellow patches indicates the state of two stable states (black dot) and one unstable state (white dot). Exchange of stability occurs at theasddel node bifurcation point (half black and half white dot) and is called the static switching point.

## 3. Theoretical formulation

Though our analysis is valid for any kind of arches with any boundary conditions, we only consider the bistable arches with fixed-fixed boundary conditions for the analysis (see Fig. xx).We consider an arch of the span $L$, the in-plane thickness $w$, out-of-plane width $d$, and made up of material which Young's modulus is $E$, and density is $\rho$. The mid-span height of the arch be taken as $h_{mid}$, $I = \frac{dw^3}{12}$ be the second moment of area, and $A = dw$ is the cross-sectional area of the profile. We assume the arch does not deform out-of-plane when we apply a force. We follow the Rayleigh-Ritz approximation to express the as-fabricated shape, $h(x)$, and the deformed shape, $w(x)$, of the arch as a weighted linear combination of the buckling mode shape of a straight column of the same length and boundary conditions. Let $\phi_i$ be the $i^{th}$ buckling mode shape, then the as-fabricated shape and the deformed shape are expressed as,

$$h(x) = h_{mid} \sum_{i=1}^{N} a_i \phi_i(x)$$
$$w(x) = h_{mid} \sum_{i=1}^{N} A_i \phi_i(x)$$
(1)

where $\hat{a}_i$ and $\hat{A}_i$ are the mode weights of the $i^{th}$ mode shape in initial, and deformed shape, respectively. Whereas $N$ is the total number of mode weights we consider to approximate the as-fabricated and deformed shape of the arch. For an fixed-fixed arch the mode weights are given by

$$\phi_i(x) = \begin{cases} 1 - \cos(M_i x/L); & i = 1, 3, 5, \ldots \\ 1 - 2x/L - \cos(M_i x/L) + 2\frac{\sin(M_i x/L)}{M_i}; & i = 2, 4, 6, \ldots \end{cases} \qquad M_i = \begin{cases} (i+1)\pi; & i = 1, 3, 5, \ldots \\ 2.86\pi, 4.92\pi, 6.94\pi; & i = 2, 4, 6, \ldots \end{cases}$$
(2)

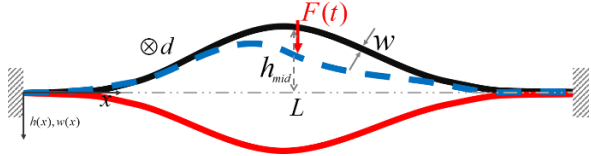

**Fig. 3.** A fixed-fixed arch with all the geometric parameters and a point load at the center in three different configurations; as-fabricated (black curve), intermediate (blue dashed curve), and inverted position (red curve) when the load exceeds the static switching value.

We take the energy method approach to derive the governing equation of the arch. Assuming the arch to be shallow (slope, $\frac{dw}{dx} \ll 1$), the bending energy of the arch is expressed as

$$se_b = \frac{EI}{2} \int_0^L \left( \frac{d^2h}{dx^2} - \frac{d^2w}{dx^2} \right)^2 dx = \frac{EI h_{mid}^2}{4L^3} \left( \sum_{i=1}^{N} A_i^2 M_i^4 + \sum_{i=1}^{N} a_i^2 M_i^4 + 2\sum_{i=1}^{N} A_i a_i M_i^4 \right)$$
(3)

When the arch deforms, it leads to axial compression. If $s_{initial}$ and $s$ be the initial and deformed arc length of the arch, respectively, the compression energy is given by

$$se_c = \frac{EA}{2L}(s_{initial} - s) = \frac{EA}{8L} \left[ \int_0^L \left( \frac{dh}{dx} \right)^2 dx - \int_0^L \left( \frac{dw}{dx} \right)^2 dx \right]^2$$
$$= \frac{EA h_{mid}^4}{32L^3} \left( \sum_{i=1}^{N} a_i^2 M_i^2 - \sum_{i=1}^{N} A_i^2 M_i^2 \right)^2$$
(4)

The work potential due to the point load $f$ at the mid-span of the is

$$wp = -f\left(1+\hat{v}t\right)\left(h\big|_{x=\frac{L}{2}} - w\big|_{x=\frac{L}{2}}\right) = -2h_{mid}f\left(1+\hat{v}t\right)\left(\sum_{i=1,5,9,\ldots} a_i - \sum_{i=1,5,9,\ldots} A_i\right)$$

(5)

Considering the inertia effect, the kinetic energy, $k_e$, of the arch is expressed as

$$k_e = \frac{1}{2}\int_0^L \rho A\left(\frac{dA_i}{dt}\right)^2 dx = \frac{\rho A L h_{mid}^2}{4}\sum_{i=1}^N \left(\frac{dA_i}{dt}\right)^2 \qquad (6)$$

here, the dot represents the derivative w.r.t. time, $t$. Assuming the arch has viscous damping, we consider the Rayleigh dissipation function

$$D = \frac{1}{2}h_{mid}^2 \hat{c}\sum_{i=1}^N \left(\frac{dA_i}{dt}\right)^2 \qquad (7)$$

where $\hat{c}$ is the damping constant of the arch. Using modified Hamilton's principle we get the governing differential equation of the arch as

$$\frac{d}{dt}\left(\frac{\partial L}{\partial\left(\partial A_i/\partial t\right)}\right) + \frac{\partial D}{\partial\left(\partial A_i/\partial t\right)} - \frac{\partial L}{\partial A_i} - Q_i^{nc} = 0 \qquad (8)$$

where $L = K_e - se_b - se_c$ is the Lagrange function and $Q_i^{nc}$ is the applied conservative and/or nonconservative forces (generalized force). Substituting the energy terms Eq. (3), (4), (6), and dissipation term from Eq. (7), in Eq. (8) we get a system of differential algebraic equations (DAEs) as

$$\frac{\rho A L h_{mid}^2}{2}\left(\frac{d^2 A_i}{dt^2}\right) + h_{mid}^2 \hat{c}\left(\frac{dA_i}{dt}\right)_i + \left(\frac{EIh_{mid}^2}{2L^3}\right)\left(A_i M_i^4 - a_i M_i^4\right) - \frac{EAh_{mid}^4}{8L^3}\left(\sum_{i=1}^N a_i^2 M_i^2 - \sum_{i=1}^N A_i^2 M_i^2\right)A_i M_i^2 - 2h_{mid}f\left(1+\hat{v}t\right)\big|_{1,5,\ldots} = 0$$

(10)

To make the analysis independent of scaling and geometric parameters, we non-dimensionalize the governing differential equations. Each term on the left hand-side of Eq. 10 represents a forcing term, and we divide each term by $\frac{EIh_{mid}}{L^3}$ to make them normalize. We introduces a dimensionless geometric parameter, $Q = \frac{h_{mid}}{w}$, and get the nondimensional goverening DAEs as Eq. (11). The other nondimensional parameters used in this paper are given in the Table. 1.

$$\frac{1}{2}\ddot{A}_i + c\dot{A}_i + \frac{1}{2}\left(A_i M_i^4 - a_i M_i^4\right) - \frac{3Q^2}{2}\left(\sum_{i=1}^N a_i^2 M_i^2 - \sum_{i=1}^N A_i^2 M_i^2\right)A_i M_i^2 - 2F\left(1+v\tau\right)\big|_{i=1,5,\ldots} = 0$$

(11)

Table. 1. Different nondimensional parameters used for the analysis of the arch.

| Physical parameter | Normalized factor | Nondimensional parameter |
|---|---|---|
| Force ($f$) | $\frac{EIh_{mid}}{L^3}$ | $F$ |
| Mode weights ($\hat{a}_i, \hat{A}_i$) | $h_{mid}$ | $a_i, A_i$ |

| Time ($t$) | $\sqrt{\dfrac{\rho A L^4}{EI}}$ | $\tau$ |
|---|---|---|
| Force rate ($\hat{\alpha}$) | $\dfrac{(EI)^{3/2} h_{mid}}{\sqrt{\rho A L^5}}$ | $\alpha$ |
| Damping co-efficient ($\hat{c}$) | $\sqrt{\dfrac{EI\rho A}{L^2}}$ | $c$ |
| | $\dfrac{h_{mid}}{w}$ | $Q$ |

## 4. Analysis of switching time by taking one-mode approximation

To understand the physics of critical slowing down of an arch, we first analyse the arch by taking a one-mode approximation for the deformed shape and then extend it to a higher mode approximation. Therefore, by taking $A_1$ in the deformed shape, we get the governing DAE of the arch as

$$\frac{1}{2}\ddot{A}_1 + c\dot{A}_1 + \frac{1}{2}\left(A_1 M_1^4 - a_1 M_1^4\right) - \frac{3Q^2}{2}\left(a_1^2 M_1^2 - A_1^2 M_1^2\right)A_1 M_1^2 - 2F(1+v\tau) = 0 \quad (12)$$

By substituting $\delta = 2(a_1 - A_1)$ in Eq. (12), we get equations in terms of $\delta$ as,

$$\ddot{\delta} + 2c\dot{\delta} + \left(3Q^2 M_1^4\right)\delta\left[\tfrac{1}{4}\delta^2 - \tfrac{3a_1}{2}\delta + \left(2a_1^2 + \tfrac{1}{3Q^2}\right)\right] - 8F(\tau) = 0 \quad (13)$$

Here, $\delta$ signifies the mid-point displacement of the arch. Eq. (13) can generalized to Eq. (14) as follows,

$$\ddot{\delta} + c\dot{\delta} = g\left(\delta, F(\tau)\right) = F(\tau) - X(\delta) \quad (14)$$

where $X(\delta) = \left(3Q^2 M_1^4\right)\left[\tfrac{1}{4}\delta^3 - \tfrac{3a_1}{2}\delta^2 + \left(2a_1^2 + \tfrac{1}{3Q^2}\right)\delta\right]$ is the internal force of the arch, and a cubic polynomial in displacement, $\delta$. Here, $F(\tau)$ is the controlling parameter that decide the stability of the arch, and vary either quasi-statically or linearly with time. When $F(\tau)$ is independent of time, there exists two stable state for $F(\tau) < F_c$, and a single stable state when $F(\tau) > F_c$. On the other hand, when $F(\tau)$ is dependent of time, there exists two stable state when $\tau < \tau_c$ such that $F(\tau) < F_c$, and a single stable state when $\tau > \tau_c$ and $F(\tau) > F_c$ such that at $F(\tau_c) = F_c$. For both the cases, the critical slowing down occurs at the staic switching point where the saddle-node bifurcation occurs and following conditions are satisfied

$$g(\delta_c, F_c) = 0, g'_\delta(\delta_c, F_c) = 0, g'_F(\delta_c, F_c) \neq 0, g''_\delta(\delta_c, F_c) \neq 0 \quad (15)$$

where $\delta_c$ and $F_c$ are the static switching displacement and force, respectively. When the controlling parameter $F$ is time-dependent, the switching point also changes with time, and the arch never switches at the static switching load. Therefore, when the load is time varying, the switching of the arch delays depending on the applied load rate. In the upcoming sections we present the analysis method, and analytical formulae for the delay switching of the arch for both time-independent, and time-dependent loading.

### 4.1 Time-independent loading

When the applied load is independent of time, Eq. (14) becomes a autonomous system and the governing equation is given as

$$\ddot{\delta} + c\dot{\delta} = g(\delta, F) = F - X(\delta) \tag{16}$$

Here the stability of the arch depends on the applied force $F$. As the force increases, the arch becomes more unstable, and completely lose the stability at the switching point $(\delta_c, F_c)$ when its saddle node bifurcation occurs. The static switching point is given as

$$\delta_c = 2a_1 - 2\sqrt{\frac{a_1^2}{3} - \frac{1}{9Q^2}}$$
$$F_c = (3Q^2 M_1^4)\left[\frac{1}{4}\delta_c^3 - \frac{3a_1}{2}\delta_c^2 + \left(2a_1^2 + \frac{1}{3Q^2}\right)\delta_c\right] \tag{17}$$

When the load is slightly more than the switching load, $F_c$, the arch switch to a remote stable state. This switching phenomena occurs at a smaller time scale compared to the elastic time scale of the arch. This slow switching is due to the faltness of the energy-displacement curve near to the static switching point, hence is a local phenomena. Therefore, to characterize the switching dynamics of the arch, we do stability analysis by linearizing Eq. (16) at the switching point $(\delta_c, F_c)$. We perturb force as $F_c + \varepsilon$ and displacement as $\delta_c + \bar{\delta}$, and use Eq. (15) to linearize Eq. (16) as follows

$$\ddot{\bar{\delta}} + c\dot{\bar{\delta}} = g(\delta_c + \bar{\delta}, F_c + \varepsilon)$$
$$\Rightarrow \ddot{\bar{\delta}} + c\dot{\bar{\delta}} = \cancel{g(\delta_c, F_c)}^0 + \cancel{g'_\delta(\delta_c, F_c)}^0 \bar{\delta} + g'_F(\delta_c, F_c)\varepsilon + g''_\delta(\delta_c, F_c)\bar{\delta}^2 + \cancel{g''_F(\delta_c, F_c)}^0 \varepsilon^2 + O(\bar{\delta}^3, \varepsilon^3)$$
$$\Rightarrow \ddot{\bar{\delta}} + c\dot{\bar{\delta}} = \varepsilon + K\bar{\delta}^2 \tag{18}$$

where $K = \left(\frac{9}{2}Q^2 M_1^4\right)\left(\frac{\delta_c}{2} - a_1\right)$ is the local curvature of the force-displacement curve at the switching point. When we consider only the first order terms, the linearization results a trivial solution. Therefore, we take the second order terms to capture the slow transition of the arch. We may consider both the inertia and damping effect simultaneously as in Eq. (18), then the analysis becomes more complicated and may not be helpful in understanding the physics of the problem. So, we analyze the switching dynamics by considering the dessipative term and the undamped term separately.

#### 4.1.1 Damping Dominated Analysis

When an arch is made up of high viscous material, the dissipation dominates over the inertial force. Therefore, the system behave as a overdamped system, and the dynamic equation of motion at the switching point for this damping dominated system is given by:

$$c\frac{d\bar{\delta}}{d\tau} = \varepsilon + K\bar{\delta}^2 \tag{19}$$

From Eq. (25) we get the relation between time, $\tau$, arch deformation parameter, $\bar{\delta}$, and force perturbation, $\varepsilon$, as follows

$$\tau = c\left(\frac{1}{\sqrt{K\varepsilon}}\right)\tan^{-1}\left(\sqrt{\frac{K}{\varepsilon}}\bar{\delta}\right) \tag{20}$$

$$\bar{\delta} = \sqrt{\frac{\varepsilon}{K}}\tan\left(\frac{\sqrt{K\varepsilon}\,\tau}{c}\right) \tag{21}$$

In the above formulations (Eq. (20) and (21)), $\varepsilon$ signifies the closeness of the applied load to the static critical load. We can observ the switching time depends on the parameter $\varepsilon$ by an inverse power law of degree $\frac{1}{2}$. Therefore, when the applied load is close to the switching load value, the switching of the arch slows down drastically. When the arch switches from one stable path to another stable path, the arch travels away from the vicinity of the switching displacement. Hence, the shape parameter, $\bar{\delta}$, blows up drastically to infinity. Therefore, for an load increment of $\varepsilon$, the arch takes a time $\tau_\infty$ to travel from static switching point to another stable point, and is given by Eq. (22).

$$\tau_\infty = c\left(\frac{1}{\sqrt{K\varepsilon}}\right)\frac{\pi}{2} = 1.5708\,c\left(\frac{1}{\sqrt{K\varepsilon}}\right) \tag{22}$$

On the other hand, the time evolution of the deformation of the arch shape after the static switching point is given by Eq. (23).

$$\bar{\delta} = \sqrt{\frac{\varepsilon}{K}}\tan\left(\frac{\sqrt{\varepsilon K}}{c}\left(\tau - \frac{\pi}{2}\right)\right) \tag{23}$$

### 4.1.2 Undamped Switching Analysis

For a non or weakly dissipative system, the damping term is negligible compered to the inertial term. The linearized differential equation for such a system to study the local behavior is given as,

$$\frac{d^2\bar{\delta}}{d\tau^2} = \varepsilon + K\bar{\delta}^2 \tag{24}$$

with initial conditions $\bar{\delta} = \frac{d\bar{\delta}}{dt} = 0$ at $\tau = 0$ at the static wsitching point. We multiply both sides of Eq. (24) with $\dot{\bar{\delta}}$, integrate twice over the time, $\tau$, and use the elliptic integral of first kind [31] to get the relation between time and the shape parameter $\bar{\delta}$ as follows

$$\int_0^\tau d\tau = \tau = \sqrt{\frac{3}{2}}\frac{1}{\sqrt{K}}\int_0^{\bar{\delta}}\frac{d\bar{\delta}}{\sqrt{\bar{\delta}\left(\bar{\delta}^2 + \frac{3\varepsilon}{K}\right)}}$$

$$\Rightarrow \tau = \frac{1}{\sqrt{2}}\left(\frac{3}{K\varepsilon}\right)^{1/4}\mathbf{F}\left(\cos^{-1}\left(\frac{\sqrt{3\varepsilon} - \sqrt{K}\bar{\delta}}{\sqrt{3\varepsilon} + \sqrt{K}\bar{\delta}}\right), \frac{1}{\sqrt{2}}\right) \tag{25}$$

where $\mathbf{F}(\bullet,\bullet)$ is the expression for elliptic intefral of first kind. We find the switching time of the arch from the static switching point to another stable state by substituting $\bar{\delta} \to \infty$ in Eq. (25) like we did in the last section. So, the critical switching time for an inertial system is given by

$$\tau_\infty = \frac{1}{\sqrt{2}} \left( \frac{3}{K\varepsilon} \right)^{1/4} \mathbf{F}\left( \pi, \frac{1}{\sqrt{2}} \right) = 3.4508 (K\varepsilon)^{-1/4} \tag{26}$$

Like the damping dominated system, the switching time of the arch obeys an inverse power law with the force perturbation, $\varepsilon$, but with a power of $\frac{1}{4}$. For both the cases, the switching of the arch delays drastically when the applied force is close to the static switching force. But, the switching of an inertia dominated system occurs at a slower time scale than a dissipative system. The switching time also depends on the local curvature, $K$, of the force-displacement curve at the switching point. Therefore, for a same force perturbation, the switching time vary for depending on its shape and boundary conditions.

### 4.2 Time-dependent loading

In section 4.1 we found out the switching time of an arch when a quasi-static load is applied. In that case, though the force chages with time, it is not significant. Therefore, it has a fixed asymptotes at the static switching point. On the other hand, when the applied load rate is significant, Eq. (14) becomes non-autonomous (Eq. 27), and has a moving asymptote. Hence, the applied load rate delay the switching and arch does not switch at the static switching load. To understand this delay switching, we take the controlling parameter, $F$, to be a linear function time. So, the governing equation (Eq. (14)) is written as

$$\ddot{\delta} + c\dot{\delta} = g(\delta, F(\tau)) = F_0 + v\tau - X(\delta) \tag{27}$$

where $F_0$ is the initial static load (if any) on the arch, $v$ is the rate at which the load is varying with time, $\tau$, and $X(\delta)$ is the internal force of the arch (defined in section 4). The time at which the applied load attain the static switching load, $F_c$, is called the critical transition time, $\tau_c$, and is expressed as $\tau_c = \frac{F_c - F_0}{v}$. To study the local behavior at the switching point, we linearize Eq. (27) by perturbing displacement, $\delta$, by an amount $\bar{\delta}$. Using Eq. (15), we get the linearized governing equation (Eq. 28) for the further study as

$$\ddot{\bar{\delta}} + c\dot{\bar{\delta}} = F_0 + v\tau - X(\delta_c) + \cancel{X'_\delta(\delta_c)\bar{\delta}}^{0} + X''_\delta(\delta_c)\bar{\delta}^2 + O(\bar{\delta}^3)$$

$$\Rightarrow \ddot{\bar{\delta}} + c\dot{\bar{\delta}} = F_0 + v\tau - F_c + K\bar{\delta}^2$$

$$\Rightarrow \ddot{\bar{\delta}} + c\dot{\bar{\delta}} = v(\tau - \tau_c) + K\bar{\delta}^2 \tag{28}$$

where $K$ is the local curvature of the force-displacement curve at the switching point. Define a dummy time scale, $\bar{\tau} = \tau - \tau_c$, and substitute in Eq. (28) to get

$$\ddot{\bar{\delta}} + c\dot{\bar{\delta}} = v\bar{\tau} + K\bar{\delta}^2 \tag{29}$$

We study the dissipative system and the undamped system separately like we did in the previous section.

### 4.2.1 Damping Dominated Analysis

When the dissipative term dominates than the inertial term, the goverening differential equation is given as

$$c\frac{d\bar{\delta}}{d\bar{\tau}} = v\bar{\tau} + K\bar{\delta}^2 \tag{30}$$

This first-order nonlinear differential equation (Eq. (30)) is called the Riccati differential equation. To solve Eq. (30) we substitute $\bar{\delta} = \frac{-c}{Ku}\left(\frac{du}{d\bar{\tau}}\right)$ and convert it to a second-order Airy differential equation (Eq. 31).

$$c^2 \frac{d^2u}{d\bar{\tau}} + (Kv)\bar{\tau}u = 0 \tag{31}$$

with initial conditions as $u(\bar{\tau}_0) = 1;\ u'(\bar{\tau}_0) = -K\bar{\delta}_0$. Here, $\bar{\tau}_0$ corresponds to the time at which we start applying the load, and $\bar{\delta}_0$ is the initial shape of the ach. The analytical solution of Eq. (31) using Airy functions is given as [33]

$$u = \left[C_1\ \mathbf{Ai}(-z) + C_2\ \mathbf{Bi}(-z)\right]\Big|_{z=\left(\frac{Kv}{c^2}\right)^{1/3}\bar{\tau}} \tag{32}$$

and its derivative w.r.t. dummy time, $\bar{\tau}$, is written as

$$u' = \frac{du}{dz}\frac{dz}{d\bar{\tau}} = -\left(\frac{Kv}{c^2}\right)^{1/3}\left[C_1\ \mathbf{Ai}'(-z) + C_2\ \mathbf{Bi}'(-z)\right]\Big|_{z=\left(\frac{Kv}{c^2}\right)^{1/3}\bar{\tau}} \tag{33}$$

where $\mathbf{Ai}(z)$, and $\mathbf{Bi}(z)$ are the Airy functions of first, and second kind, respectively. $C_1$, and $C_2$ are the movable parameteric constant and defined in terms of the initial condition as

$$\begin{aligned} C_1 &= \frac{(K)^{2/3} v^{-1/3}\bar{\delta}_0\ \mathbf{Bi}(z_0) - \mathbf{Bi}'(z_0)}{\mathbf{Ai}(z_0)\mathbf{Bi}'(z_0) - \mathbf{Ai}'(z_0)\mathbf{Bi}(z_0)}\Big|_{z_0=\left(\frac{Kv}{c^2}\right)^{1/3}\bar{\tau}_0} \\ C_2 &= \frac{(K)^{2/3} v^{-1/3}\bar{\delta}_0\ \mathbf{Ai}(z_0) - \mathbf{Ai}'(z_0)}{\mathbf{Ai}(z_0)\mathbf{Bi}'(z_0) - \mathbf{Ai}'(z_0)\mathbf{Bi}(z_0)}\Big|_{z_0=\left(\frac{Kv}{c^2}\right)^{1/3}\bar{\tau}_0} \end{aligned} \tag{34}$$

Using Eq. (32) and Eq. (33) we solve Eq. (30) to get the arch deformation parameter as

$$\bar{\delta} = \frac{-c}{Ku}\left(\frac{du}{d\bar{\tau}}\right) = (K)^{-2/3} c^{1/3} v^{1/3}\left[\frac{C_1\ \mathbf{Ai}'(z) + C_2\ \mathbf{Bi}'(z)}{C_1\ \mathbf{Ai}(z) + C_2\ \mathbf{Bi}(z)}\right]\Big|_{z=\left(\frac{Kv}{c^2}\right)^{1/3}\bar{\tau}}$$

$$\Rightarrow \bar{\delta} = (K)^{-2/3} c^{1/3} v^{1/3}\left[\frac{\mathbf{Ai}'(z) + \frac{C_2}{C_1}\mathbf{Bi}'(z)}{\mathbf{Ai}(z) + \frac{C_2}{C_1}\mathbf{Bi}(z)}\right]\Big|_{z=\left(\frac{Kv}{c^2}\right)^{1/3}\bar{\tau}} \tag{35}$$

In equation Eq. (35) the evolution of the deformed shape depends on the ratio $\frac{C_1}{C_2}$, and hence on the initial conditions of the arch. As we are considering the large-time behavior of a non-autonomous dynamic system, instead of going forward in the time, $\bar{\tau} \to +\infty$, we rewind time such that the initial state corresponds to the time $\bar{\tau} \to -\infty$ (see Fig. (4)). Here, the initial state of the arch act as a pullback attractor of the system.

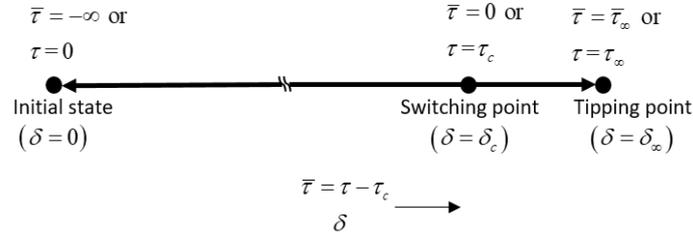

**Fig. 4.** Concept of a pullback attractor. The initial state of the arch act as a pullback attractor w.r.t. the static switching point of the arch and can be assumed as a state that achieved long back at $\bar{\tau} \to -\infty$ before reaching the switching point.

As the initial state of the arch corresponds to $\bar{\tau}_0 \to -\infty$, $z \to -\infty$, hence $\mathbf{Ai}(z) \to 0$, $\mathbf{Bi}(z) \to \infty$, and the ratio $\frac{C_1}{C_2} \to 0$ (Eq. 36).

$$\lim_{z_0 \to -\infty}\left(\frac{C_2}{C_1}\right) = \lim_{z_0 \to -\infty}\left(\frac{(K)^{2/3} v^{-1/3}\bar{\delta}_0 \mathbf{Ai}(z_0) - \mathbf{Ai}'(z_0)}{(K)^{2/3} v^{-1/3}\bar{\delta}_0 \mathbf{Bi}(z_0) - \mathbf{Bi}'(z_0)}\right)\Bigg|_{z_0 = \left(\frac{Kv}{c^2}\right)^{1/3}\bar{\tau}_0} = 0 \tag{36}$$

Using Eq. (36) in Eq. (35), we get the expression for $\bar{\delta}$ reduces to

$$\bar{\delta} = (K)^{-2/3} c^{1/3} v^{1/3} \left[\frac{\mathbf{Ai}'(-z)}{\mathbf{Ai}(-z)}\right]_{z=\left(\frac{Kv}{c^2}\right)^{1/3}\bar{\tau}} \tag{37}$$

and the expression for $\delta$ is given as

$$\delta = \delta_c + (K)^{-2/3} c^{1/3} v^{1/3} \left[\frac{\mathbf{Ai}'(-z)}{\mathbf{Ai}(-z)}\right]_{z=\left(\frac{Kv}{c^2}\right)^{1/3}\bar{\tau}} \tag{38}$$

When the arch switches from one stable path to the other the deformation parameter becomes large ($\bar{\delta} \to \infty$). As we can observe, the evolution of the arch shape depends on the function $\frac{\mathbf{Ai}'(-z)}{\mathbf{Ai}(-z)}$ (Eq. (37) and (38)). This function shows an asymptotic behavior as $z \to 2.3381$. At this point $\mathbf{Ai}(-z) \to 0$ hence the function $\frac{\mathbf{Ai}'(-z)}{\mathbf{Ai}(-z)} \to \infty$ (see Fig. xx). Therefore, the switching time of the arch, $\bar{\tau}_\infty$, corresponds to $z = 2.3381$, and expressed in Eq. (39).

$$z = \left(\frac{Kv}{c^2}\right)^{1/3} \bar{\tau}_\infty = 2.3381$$

$$\Rightarrow \bar{\tau}_\infty = 2.3381 \left(\frac{Kv}{c^2}\right)^{-1/3} \tag{39}$$

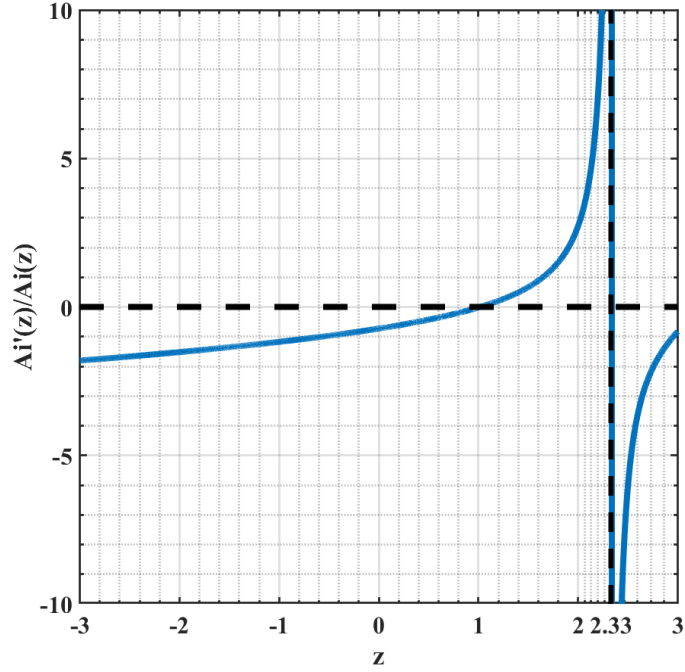

**Fig. 5.** Shows the behavior of the ratio $\frac{Ai'(z)}{Ai(z)}$ on the real axis. It shows a singularity at $z = 2.3381$.

The time an arch takes to switch from the static switching point to the tipping point is defined as $\bar{\tau}_\infty$. We get the tipping point time, $\tau_\infty$, of the arch by changing the dummy time scale, $\bar{\tau}$, back to actual time scale, $\tau$.

$$\tau_\infty = \tau_c + 2.3381 \left( \frac{Kv}{c^2} \right)^{-\frac{1}{3}} \tag{40}$$

From Eq. (40) we can see, the arch does not switch to another stable point when the force reach the staic switching value. Instead, the arch takes some additionl load, and the switching delay depending on the load rate and damping of the system. The amount of extra load an arch takes is called the delay of the system.

### 4.2.2 Undamped Switching Analysis

For a linearly varying load, the linearized governing differential equation for a non-dissipative system is given as

$$\frac{d^2 \bar{\delta}}{d \bar{\tau}^2} = v(\bar{\tau}) + K \bar{\delta}^2 \tag{41}$$

This equation (Eq. (41)) is a painleve equation of first kind. This equation cannot be solved to get a general solution for any arbitrary time. But it is asymptotic for a large negative time, i.e., for $\bar{\tau} \to -\infty$. For a large negative time domain, the solution of Eq. (41) is oscillatory about the static equilibrium point. We simplify Eq. (41) by changing time, $\bar{\tau}$, and defrmation, $\bar{\delta}$, to Boutroux's coordinates [32] as

$$\bar{\delta} = \sqrt{\bar{\tau}} Y$$
$$X = \frac{4}{5}(\bar{\tau})^{5/4} \tag{42}$$

where $X$ and $Y$ are the dummy dependent and independent variables. Substituting Eq. (42) in Eq. (41), we get the simplified governing equation as

$$\frac{d^2Y}{dX^2} + \frac{1}{X}\frac{dY}{dX} - \frac{4}{25}\frac{1}{X^2}Y = v + KY^2 \tag{43}$$

Like the previous analysis, we focous on the large-time behaviour of the arch the concept of pullback attractor to solve it. So, for large negative time, i.e., for $\bar{\tau} \to -\infty$, the relation $\lim_{\bar{\tau} \to -\infty}\left(\frac{1}{X}\right) \to 0$, and $\lim_{\bar{\tau} \to -\infty}\left(\frac{1}{X^2}\right) \to 0$ holds. Therefore, Eq. (43) reduce to

$$\frac{d^2Y}{dX^2} = v + KY^2 \tag{44}$$

with initial conditions $Y = \frac{dY}{dX} = 0$ at $X = 0$. We use elliptic integral of first kind to get the relation between the $X$ and $Y$ as follows

$$\int_0^X dX = X = \sqrt{\frac{3}{2}}\frac{1}{\sqrt{K}}\int_0^Y \frac{dY}{\sqrt{Y\left(Y^2 + \frac{3v}{K}\right)}}$$

$$\Rightarrow X = \frac{1}{\sqrt{2}}\left(\frac{3}{K}\right)^{1/4}\left(v^{-\frac{1}{4}}\right)F\left(\cos^{-1}\left(\frac{\sqrt{3v} - \sqrt{KY}}{\sqrt{3v} + \sqrt{KY}}\right), \frac{1}{\sqrt{2}}\right) \tag{45}$$

When the arch switch form the static switching point to another stable point, the deformation parameter, $\bar{\delta}$, becomes very large hence $Y \to \infty$. So, Eq. (45) gives

$$X_\infty = 3.4508\left(\frac{1}{Kv}\right)^{1/4} \tag{46}$$

By transferring Boutroux's coordinate to the original coordinate using Eq. (42), we get

$$\bar{\tau}_\infty = 3.22\left(\frac{1}{Kv}\right)^{1/5} \tag{47}$$

and the total time to switch from undeformed state to the tipping point is given by

$$\tau_\infty = \tau_c + 3.22\left(\frac{1}{Kv}\right)^{1/5} \tag{48}$$

As we can see, the delay in the switching of the arch obeys a power law with the applied load rate, $v$, and with the local curvature, $K$, of the arch. For an particular arch, the amount of delay increases with a increase in the applied load rate. So, when the load is time-dependent, the arch never switches at the staic switching load even for a slowly varying load.

## 5. Analytical solution of switching time taking multi-mode approximation

In Section 4 we took only the first mode to approximate the deformed shape and predict the delay switching of the arch. Taking the higher modes to approximate the deformed shape, improve the accuracy of the analysis. When we consider higher modes, we solve a system of coupled nonlinear ADEs to predict the delay switching of the arch. Though we can take any number of modes in the approximation, in this section we take a two-mode-approximation model (see Eq. (49)) for the analysis.

The ADEs of a two-mode-approximation model for a time-independent loading and damping dominated system by taking first and fifth mode weight in the as-fabricated and deformed shape of the arch is expressed as

$$c \frac{dA_1}{d\tau} = -2F - f(A_1, A_5)$$
$$c \frac{dA_5}{d\tau} = -2F - g(A_1, A_5)$$
(49)

where $f$ and $g$ signifies the internal forces of the arch and is given by

$$f(A_1, A_5) = \frac{1}{2}(A_1 M_1^4 - a_1 M_1^4) - \frac{3Q^2}{2}(a_1^2 M_1^2 + a_5^2 M_5^2 - A_1^2 M_1^2 - A_5^2 M_5^2)(A_1 M_1^2)$$
$$g(A_1, A_5) = \frac{1}{2}(A_5 M_5^4 - a_5 M_5^4) - \frac{3Q^2}{2}(a_1^2 M_1^2 + a_5^2 M_5^2 - A_1^2 M_1^2 - A_5^2 M_5^2)(A_5 M_5^2)$$
(50)

We get the switching mode weights ($A_{1c}, A_{5c}$) and the switching load value ($F_c$) by finding the peak in the force-displacement curve of the arch [34]. To study the delay switching of the arch, we linearize Eq. (50) around the critical switching point by taking $A_1 \to A_{1c} + \overline{A}_1$ and $A_5 \to A_{5c} + \overline{A}_5$ to get the expressions as

$$f(A_{1c} + \overline{A}_1, A_{5c} + \overline{A}_5) = -2F_c + \nabla \mathbf{f}\big|_{(A_{1c},A_{5c})} \begin{Bmatrix} \overline{A}_1 \\ \overline{A}_5 \end{Bmatrix} + \frac{1}{2}\begin{Bmatrix} \overline{A}_1 \\ \overline{A}_5 \end{Bmatrix}^T \mathbf{H}_f\big|_{(A_{1c},A_{5c})} \begin{Bmatrix} \overline{A}_1 \\ \overline{A}_5 \end{Bmatrix} + h.o.t$$
$$g(A_{1c} + \overline{A}_1, A_{5c} + \overline{A}_5) = -2F_c + \nabla \mathbf{g}\big|_{(A_{1c},A_{5c})} \begin{Bmatrix} \overline{A}_1 \\ \overline{A}_5 \end{Bmatrix} + \frac{1}{2}\begin{Bmatrix} \overline{A}_1 \\ \overline{A}_5 \end{Bmatrix}^T \mathbf{H}_g\big|_{(A_{1c},A_{5c})} \begin{Bmatrix} \overline{A}_1 \\ \overline{A}_5 \end{Bmatrix} + h.o.t$$
(51)

where $\mathbf{H}_f$ and $\mathbf{H}_g$ be the Hessian matrix of $f$ and $g$, respectively at the critical point. We also perturb the applied load by an amount $\varepsilon$, i.e., $F \to F_c + \varepsilon$. We substitute Eq. (51) in Eq. (49) to get the linearized equation

$$c \frac{d}{d\tau}\begin{Bmatrix} \overline{A}_1 \\ \overline{A}_5 \end{Bmatrix} = \begin{Bmatrix} \varepsilon \\ \varepsilon \end{Bmatrix} + \begin{bmatrix} \nabla \mathbf{f}\big|_{(A_{1c},A_{5c})} \\ \nabla \mathbf{g}\big|_{(A_{1c},A_{5c})} \end{bmatrix} \begin{Bmatrix} \overline{A}_1 \\ \overline{A}_5 \end{Bmatrix} + \begin{bmatrix} \frac{1}{2}\begin{Bmatrix} \overline{A}_1 \\ \overline{A}_5 \end{Bmatrix}^T \mathbf{H}_f\big|_{(A_{1c},A_{5c})} \begin{Bmatrix} \overline{A}_1 \\ \overline{A}_5 \end{Bmatrix} \\ \frac{1}{2}\begin{Bmatrix} \overline{A}_1 \\ \overline{A}_5 \end{Bmatrix}^T \mathbf{H}_g\big|_{(A_{1c},A_{5c})} \begin{Bmatrix} \overline{A}_1 \\ \overline{A}_5 \end{Bmatrix} \end{bmatrix}$$
(52)

As we have already dicussed, the energy landscape of the arch has a saddle point at the switching point. At this point one of the eigenvalues become zero [35]. In other words, the gradient matrix $\begin{bmatrix} \nabla \mathbf{f}\big|_{(A_{1c},A_{5c})} \\ \nabla \mathbf{g}\big|_{(A_{1c},A_{5c})} \end{bmatrix}$ has a zero eigenvalue at ($A_{1c}, A_{5c}$). If $\mathbf{V}_1$ be the unit eigenvector corresponding to the zero eigenvalue, we get the following relationship

$$\left[\begin{array}{c} \nabla \mathbf{f}\big|_{(A_{1c},A_{5c})} \\ \nabla \mathbf{g}\big|_{(A_{1c},A_{5c})} \end{array}\right] \mathbf{V}_1 = \lambda_1 \mathbf{V}_1$$

$$\Rightarrow \left\{\begin{array}{c} \bar{A}_1 \\ \bar{A}_5 \end{array}\right\}^T \left[\begin{array}{c} \nabla \mathbf{f}\big|_{(A_{1c},A_{5c})} \\ \nabla \mathbf{g}\big|_{(A_{1c},A_{5c})} \end{array}\right] \mathbf{V}_1 = \lambda_1 \left\{\begin{array}{c} \bar{A}_1 \\ \bar{A}_5 \end{array}\right\}^T \mathbf{V}_1 \quad \Rightarrow \mathbf{V}_1^T \left[\begin{array}{c} \nabla \mathbf{f}\big|_{(A_{1c},A_{5c})} \\ \nabla \mathbf{g}\big|_{(A_{1c},A_{5c})} \end{array}\right] \left\{\begin{array}{c} \bar{A}_1 \\ \bar{A}_5 \end{array}\right\} = \lambda_1^{\,0} \mathbf{V}_1^T \left\{\begin{array}{c} \bar{A}_1 \\ \bar{A}_5 \end{array}\right\}$$

$$\Rightarrow \left[\begin{array}{c} \nabla \mathbf{f}\big|_{(A_{1c},A_{5c})} \\ \nabla \mathbf{g}\big|_{(A_{1c},A_{5c})} \end{array}\right] \left\{\begin{array}{c} \bar{A}_1 \\ \bar{A}_5 \end{array}\right\} = \mathbf{0} \tag{53}$$

From Eq. (53) in Eq. (52) we get the reduced nonlinear ADEs as

$$c \frac{d}{d\tau} \left\{\begin{array}{c} \bar{A}_1 \\ \bar{A}_5 \end{array}\right\} = \left\{\begin{array}{c} \varepsilon \\ \varepsilon \end{array}\right\} + \left[\begin{array}{c} \dfrac{1}{2}\left\{\begin{array}{c} \bar{A}_1 \\ \bar{A}_5 \end{array}\right\}^T \mathbf{H}_f\big|_{(A_{1c},A_{5c})} \left\{\begin{array}{c} \bar{A}_1 \\ \bar{A}_5 \end{array}\right\} \\ \dfrac{1}{2}\left\{\begin{array}{c} \bar{A}_1 \\ \bar{A}_5 \end{array}\right\}^T \mathbf{H}_g\big|_{(A_{1c},A_{5c})} \left\{\begin{array}{c} \bar{A}_1 \\ \bar{A}_5 \end{array}\right\} \end{array}\right] \tag{54}$$

The general solution of Eq. (54) can be expanded by the combinations of the eigenvalues and eigenvectors of the gradient matrix. As the gradient matrix is singular (as one of the eigenvalues is zero), and we are looking for the long-term behavior of the system, the eigenvector corresponds to the zero eigenvalue has a greater influence on the solution of the system. Hence, the solution of Eq. (54) can be expressed as

$$\left\{\begin{array}{c} \bar{A}_1 \\ \bar{A}_5 \end{array}\right\} = \bar{\delta}(\tau) \mathbf{V}_1 \tag{55}$$

where $\bar{\delta}$ is a time-dependent scalar parameter, and is defined as the ratio of the change of mode weights to the eigenvector corresponds to the zero eigenvalue. Substituting Eq. (55) in Eq. (54) results,

$$c \frac{d\bar{\delta}}{d\tau} \mathbf{V}_1 = \left\{\begin{array}{c} \varepsilon \\ \varepsilon \end{array}\right\} + \bar{\delta}^2 \left[\begin{array}{c} \dfrac{1}{2} \mathbf{V}_1^T \mathbf{H}_f\big|_{(A_{1c},A_{5c})} \mathbf{V}_1 \\ \dfrac{1}{2} \mathbf{V}_1^T \mathbf{H}_g\big|_{(A_{1c},A_{5c})} \mathbf{V}_1 \end{array}\right] \tag{56}$$

Using Eq. (55) we convert a two variable problem in $\bar{A}_1$ and $\bar{A}_5$ to a system of ADEs in single variable, $\bar{\delta}$. We multiply both side by $\mathbf{V}_1^T$ to get an equation equivalent to the standard form as we describe in Section 4.1.1 (see Eq. (19)).

$$c \frac{d\bar{\delta}}{d\tau} = \bar{\varepsilon} + \bar{K} \bar{\delta}^2 \tag{63}$$

where $\bar{\varepsilon} = \varepsilon \operatorname{Tr}(\mathbf{I}\mathbf{V}_1)$, $\mathbf{I}$ is the identity matrix, and $\bar{K} = \mathbf{V}_1^T \left[\begin{array}{c} \dfrac{1}{2} \mathbf{V}_1^T \mathbf{H}_f\big|_{(A_{1c},A_{5c})} \mathbf{V}_1 \\ \dfrac{1}{2} \mathbf{V}_1^T \mathbf{H}_g\big|_{(A_{1c},A_{5c})} \mathbf{V}_1 \end{array}\right]$ is the

revised curvature of the force-displacement curve at the switching point. So, the expression for the switching time of the arch from the critical switching point to another stable state is expressed as

$$\tau_\infty = c\left(\frac{1}{\sqrt{\bar{K}}}\right)\left(\frac{1}{\sqrt{\bar{\varepsilon}}}\right)\frac{\pi}{2} = 1.5708\, c\left(\frac{1}{\sqrt{\bar{K}}}\right)\left(\frac{1}{\sqrt{\bar{\varepsilon}}}\right) \qquad (64)$$

We can follow the same mathematical steps for any number of modes and any loading conditions to get the critical switching time as given in the following table (see Table 2).

Table. 2. Nondimensional switching time of the arch taking higher modes in the approximation.

| | | |
|---|---|---|
| Case – 1<br>Time-independent load | Damping dominated analysis | $\tau_\infty = 1.5708\, c\left(\frac{1}{\sqrt{\bar{K}}}\right)\left(\frac{1}{\sqrt{\bar{\varepsilon}}}\right)$ |
| | Undamped analysis | $\tau_\infty = 3.4508\left(\bar{K}\bar{\varepsilon}\right)^{-1/4}$ |
| Case – 2<br>Time-dependent load (Load is Ramp function of time) | Damping dominated analysis | $\tau_\infty = \tau_c + 2.3381\left(\frac{\bar{K}v}{c^2}\right)^{-1/3}$ |
| | Undamped analysis | $\tau_\infty = \tau_c + 3.22\left(\frac{1}{\bar{K}v}\right)^{1/5}$ |

As we observe, though the scaling parameter changes as we consider higher modes for the approximation, the switching time follow the same power law with the modified local curvature and the perturbation parameters. This delay switching is an inherent property of nonlinear arches, and occurs due to the exchange of stability at the critical switching point. The detailed numerical results and the analytical results are presented in the next section.

## 6. Results and Discussion

We verified the accuracy of the analytical expressions with the numerical solutions of the respective governing equations. We numerically integrate the governing differential equations by RK4 method using the ODE45 toolbox in MATLAB. We took the tolerance and the absolute error of $10^{-12}$ for the numerical analysis. For a damping dominated system, the nondimensionalized damping, $c$, is taken to be 100 for all the cases. The displacement and velocity at the critical switching point were taken as an initial conditions for a time-independent loading. Whereas for a time-dependent loading, the undeformed shape of the arch with zero velocity and zero initial force were taken as the initial conditions.

In Fig. 4(a) and (b), we can observe for the time-independent loading, the analytical and the numerical result matches with each other within $\pm 1\%$ error for both damping and undamped analysis. It is also observed that for any choice of $Q$ value, the switching time of the arch increases with decrease in the perturbation parameter, $\varepsilon$. Therefore, the arch slows down drastically near the static switching point, and the switching time decreases as the arch moves away form the switching point. It is also observed that for a constant perturbation parameter, $\varepsilon$, the switching time decreases as the $Q$ value increases. This is because, as $Q$ increases, the value of the local stiffness at the switching point increases, hence the energy also increases. Therefore, the arch releases the store energy more rapidly by switching from one state to other, and the switching time of the arch decreases.

It is seen that for a same initial configuration, the single mode approximation underestimates the switching time for both damping dominated and undamped case when compared to the switching time by considering the higher modes. This is due to the higher local stiffness value at the static switching point when we consider one mode approximation. As we are more interested towards designing the arches, consideration of the higher modes has an advantage over the single mode while designing the arches. Fig. 4(d) shows the contour of the switching time for different initial shape and the $Q$ value for a damping dominated arch. It is observed that with increase in the magnitude of the higher mode weights in the as-fabricated shape, the switching time increases due to the decrease in the local stiffness value at the static switching point of the arch. It is noted that while designing arches for a specific switching time, the shape and geometry of the arch may not to be unique. We can get infinitely many solutions for a same switching time. Therefore, the selection of arch dimensions depends on the design constraints such as, length, weight, switching force, comeback time of the arch.

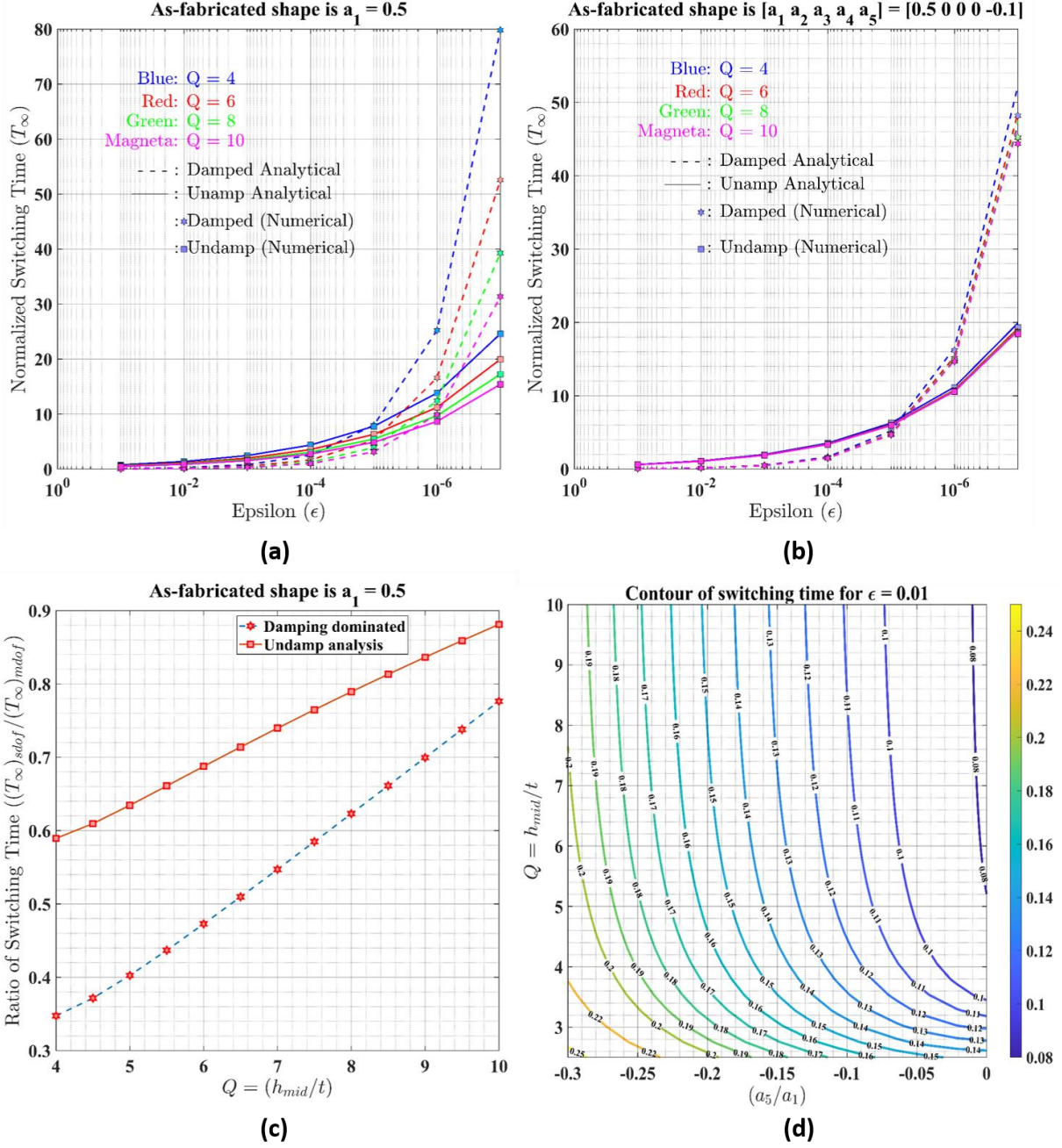

**Fig. 6.** (a) Comparision of analytical and numerical switching time of an arch for a time-independent loading for different $Q$ and perturbation parameter, $\varepsilon$, by taking one mode approximation, (b) comparision of analytical and numerical switching time by taking higher mode approximation for different $Q$ and perturbation parameter, $\varepsilon$, (c) comparision of switching time of the arch between one mode and higher mode approximation for different $Q$ value, (d) contour of switching time of the arch for different shape and different $Q$ value for $\varepsilon = 10^{-2}$. Damping constant, $c$, taken to be 100 for the analysis.

When the load is time-dependent, we do the numerical anlysis for different load rate to get some insight to the dynamic behavior of the arch. This dynamic charecteristics are compared with the static force-displacement curve and shown in Fig. 5. For both damping dominated system and undamped system, as we increase the load, the arch follows the static load-deflection curve till the vicinity of the static switching point. After this point the dynamic force-displacement curve of the starts deviating, looses its stability, and jump to the other branch of the static curve (see Fig. 5(a)-(d)), called the tipping point. In Fig. 5(b) and (d), we can observe that while switching, the load is higher than the static switching force no matter how slowly we vary the load. The switching force increases with increase in the load rate, and can approximately expressed as $F_{switch} \simeq 1.02 \times t_c \times v$. Also, the tipping point of the arch shifts farther as the load rate increases, and continues to follow the static force-deflection curve. Unlike the damped case, the load-deflection curve for undamped case undulates around the static curve at the initial (see Fig. 5(e)) and post switching (see Fig. 5(d)) period due to the absence of any dissipaton energy. As we can see, the amplitude of the undulation increases as the applied load rate increases.

    Fig. 6(a) shows the comparison of the analytical and numerical results for damped and undamped cases for time-dependent loading condition. It can be observed that the switching time for both damped and undamped case coincide with each other. The switching time decreases as the load rate increases as the arch attains its static switching load sooner for a higher load rate. Fig. 6(b) shows the comparison of the switching time for a one-mode and multi-mode approximation. Like the time-independent loading condition, the single mode approximation result underestimates the switching time of the arch. Fig. 6(c) shows the contour of switching time of the arch for different shape and applied load rate for $Q=6$. We can observe the switching time of the arch is minimum when the mode ratio, $\frac{a_5}{a_1} \simeq 0.8$, and the switching time increases as this ratio changes. A similar behavior can also be observed in the switching force contour of the arch (see Fig. 6(d)). The arch experience a maximum switching force for a mode ratio, $\frac{a_5}{a_1} \simeq 0.8$, and the switching force decreases as this ratio changes. In Fig. 6(d), we can also observe that for arch of particular shape, the switching force value increases with increase in the applied load rate. Therefore, an arch never switch at its static switching force value when a linearly varying load is applied. Fig. 6(c) and (d) can be used to design arches for different switching time and switching force requirement.

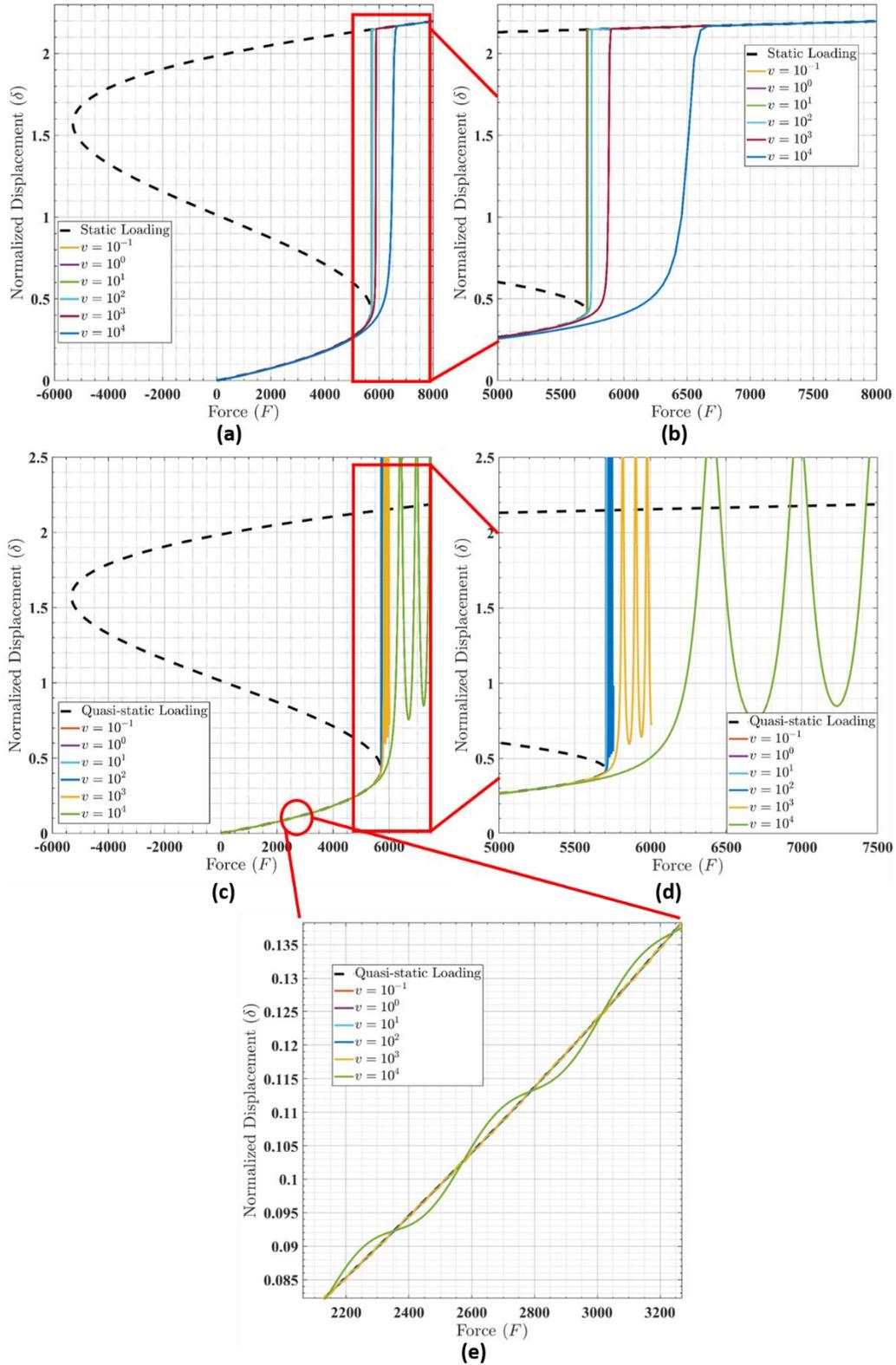

**Fig. 7.** Numerical results of damping dominated and undamped system for time-dependent loading. (a) comparision of dynamic force-displacement curve (solid) with the static curve (dash) for different load rate, $v$, of a damping dominated system (b) close-up view of the red-rectangular region showing the delay behavior for different load rate, (c) comparison of dynamic force-displacement curve (solid) with the static curve (dash) for different load rate, $v$, of a undamped system (d) close-up view of the red rectangular region showing the delay behavior for different load rate and the post switching oscilation of dynamic curve around the staic curve, (e) close-up view of the pre-switching oscillation (red circle) of the dynamic curve for different load rate. Damping constant, $c$, taken to be $100$ for the analysis.

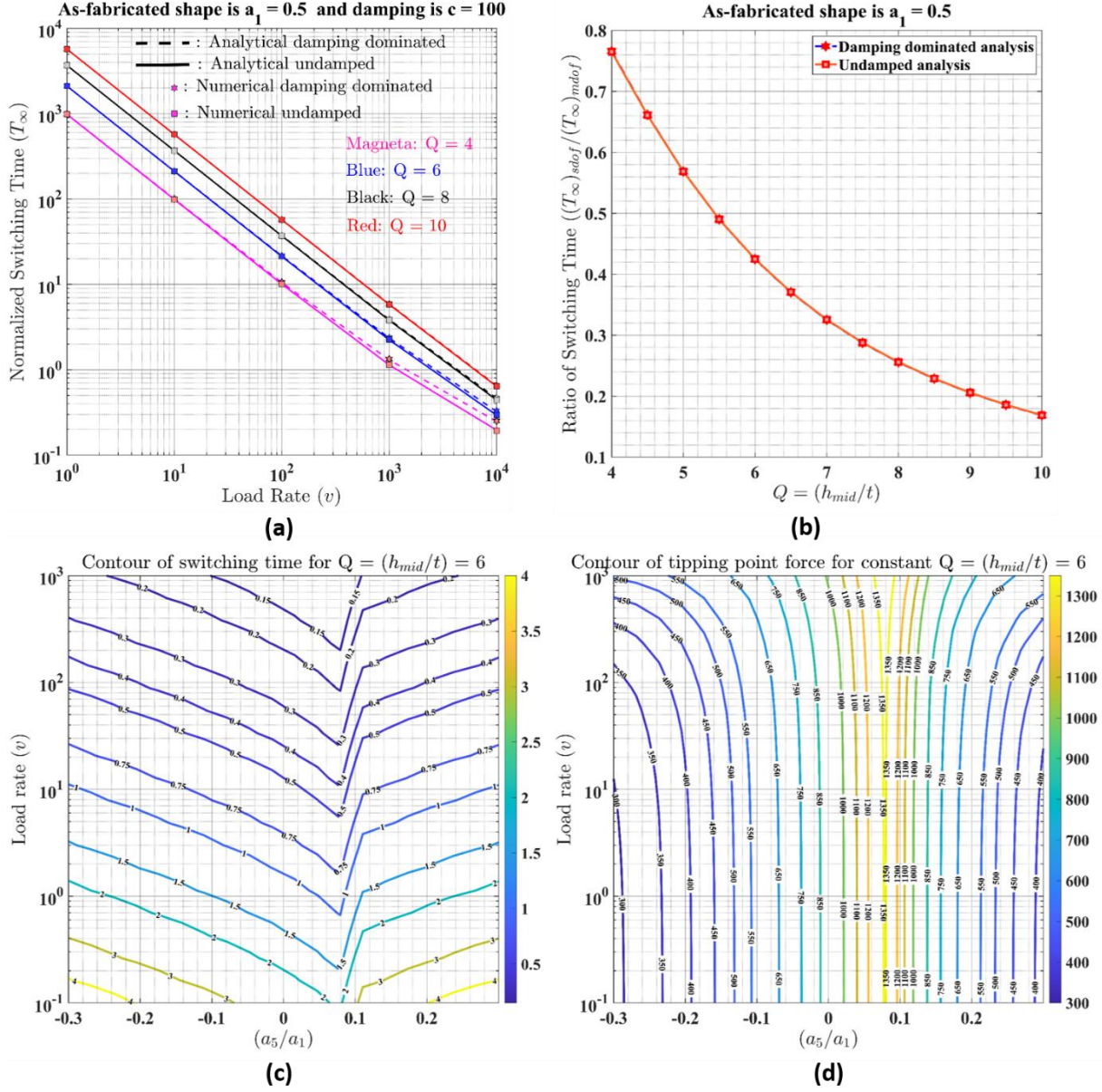

**Fig. 8.** (a) Comparision of analytical and numerical switching time of an arch for a time-dependent loading for different $Q$ and load rate, $v$, by taking one mode approximation, (b) comparison of switching time of the arch between one mode and higher mode approximation for different $Q$ value, (c) contour of switching time of the arch for different shape and different load rate for $Q=6$, (d) contour of dynamic switching force of different arch for different load rate for $Q=6$ when a ramp load is applied on the arch. Damping constant, $c$, taken to be $100$ for the analysis.

## 7. Summary

In this paper, we discussed the delay switching of shallow arches for time-dependent and time-independent loading conditions. The heighlights of the paper are as follows:

- Analytical models were developed that capture the slowing down phenomena of nonlinear shallow arches for both time-independent and linearly varying loading conditions.

- For time-independent loading, the switching time obeys an inverse power law with the local stiffness value and the force perturbation at the static switching point.
- As we get closer to the static switching point, the arch slows down drastically, and the switching time increases as the applied load gets closer to the staic switching force.
- For time-dependent loading, the switching time obeys an inverse power law with the local stiffness value as well as with the applied load rate.
- As we increase the load rate, the force taken by an arch to switch from one stable path to the other increses. Also, the tipping poind shifts farther for a increase in the applied load rate. Therefore, when the load is time varying an arch never switches at its staic switching force value.
- An analytical model was also developed to capture the delay switching of the arch when higher modes are considered to approximate the as-fabricated and deformed shapes of the arch.
- Both analytical results were compared with the numerical solutions obtained using RK4 method, and they show a fair agreement between them.
- Though we have only considered the symmetric mode of buckling in the analysis, the developed analytical method can also be used for an arch showing asymmetric bifurcation.
- The developed formulation not only valid for bistable arches, but also for the arches showing snap-through characteristic.

We can use these analytical results to design arches for different engineering applications where the applied load is time dependent. One particular application we are thinking of is in self-offloading footwear [8]. Using the analytical framework, we can coustomize the footwear based on person's weight and walking speed.

**Acknowledgements**